\documentclass[11pt, final]{amsart}
\usepackage{amssymb}
\usepackage[mathscr]{eucal}

\theoremstyle{plain}
\newtheorem{theorem}{Theorem}

\newtheorem*{main}{Main Theorem}

\theoremstyle{remark}

\def\C{\mathbb C}

\def\N{\mathbb N}
\def\ot{\otimes}
\def\R{\mathbb R}

\def\si{\sigma}
\def\supp{\operatorname{supp}}
\def\tr{\operatorname{tr}}

\begin{document}
\title[]{Some characterizations of the automorphisms of $B(H)$ and
$C(X)$} \author{LAJOS MOLN\'AR}
\address{Institute of Mathematics\\
         Lajos Kossuth University\\
         4010 Debrecen, P.O.Box 12, Hungary}
\email{molnarl@math.klte.hu}
\thanks{  This research was supported by the
          Hungarian National Foundation for Scientific Research
          (OTKA), Grant No. T030082, T031995, and by
          the Ministry of Education, Hungary, Reg.
          No. FKFP 0349/2000}
\date{\today}
\begin{abstract}
We present some nonlinear characterizations of the automorphisms of the
operator algebra $B(H)$ and the function algebra $C(X)$ by
means of their spectrum preserving properties.
\end{abstract}
\maketitle

\section{Introduction}

Surjective linear maps between Banach algebras which preserve the
spectrum are extensively studied in connection with a
longstanding open problem sometimes called Kaplansky's problem on
invertibility preserving linear maps. A weeker version of that problem
reads as follows.
Is it true that between semisimple Banach algebras every surjective
linear map which preserves the spectrum is a Jordan homomorphism?
For the algebra $B(X)$ of all bounded linear operators acting
on a Banach space this was proved to be true by Jafarian and
Sourour in \cite{JafSou}. As for commutative semisimple Banach algebras
(for instance, the algebra $C(X)$ of all continuous complex functions on
the compact Hausdorff space $X$) we once again have affirmative answer
to the question.
Namely, in that case the result is a trivial consequence of the famous
Gleason-Kahane-\. Zelazko theorem characterizing multiplicative
linear functionals.

The aim of this paper is to investigate a similar problem omitting the
condition of linearity. Clearly, nonlinear spectrum preserving
transformations can be almost arbitrary. So, we have to impose a more
restrictive condition. This will be the following: We assume that
the spectrum of the product of the images of any two elements is equal
to the spectrum of the product of that two elements. We shall see that
on the studied algebras those transformations are ''almost''
automorphisms.
Furthermore, we consider another
preserving condition concerning the spectrum which will turn out to be
in
close relation to *-automorphisms. More precisely, the main results of
the paper can be summarized as follows.

\begin{main}
Let $H$ be an infinite dimensional Hilbert space.

If $\phi: B(H) \to B(H)$ is a surjective function with the property that
\[
\si(\phi(A)\phi(B))=\si(AB) \qquad (A,B\in B(H)),
\]
then $\phi$ is either an algebra automorphism or the negative of an
algebra automorphism of $B(H)$.

If $\psi: B(H) \to B(H)$ is a surjective function with the property that
\[
\si(\psi(A)^*\psi(B))=\si(A^*B) \qquad (A,B\in B(H)),
\]
then $\phi$ is an algebra *-automorphism of $B(H)$ multiplied by a
fixed unitary element.

If $H$ is finite dimensional, then in addition to the possibilities
above we also get that $\phi$ can be an algebra antiautomorphism
or the negative of an algebra antiautomorphism of $B(H)$ and $\psi$ can
be an algebra *-antiautomorphism multiplied by a fixed unitary element.

If $X$ is a first countable compact Hausdorff space and
$\phi :C(X) \to C(X)$ is a surjective function with the property that
\[
\si(\phi(f)\phi(g))=\si(fg) \qquad (f,g \in C(X)),
\]
then $\phi$ is an algebra automorphism of $C(X)$
multiplied by a fixed continuous real function of modulus 1.

If
$\psi :C(X) \to C(X)$ is a surjective function with the property that
\[
\si({\overline{\psi(f)}}\psi(g))=\si(\overline fg) \qquad (f,g \in
C(X)),
\]
then $\psi$ is an algebra (*-)automorphism of $C(X)$ multiplied by a
fixed continuous complex function of modulus 1.
\end{main}

The statement follows from the results of the paper which follow.
We note that the referee kindly informed us about recent results on
spectrum preserving maps which are not assumed to be linear:
see \cite{BarRan} and also \cite{Hoch}.
Furthermore, we remark that other nonlinear characterizations of
the automorphisms
of matrix algebras using preserving properties can be found in
\cite{Petek}.

\section{Results}

We first fix the notation and definitions that we shall use throughout.

Every linear space is considered over the complex field.
Let $X$ be a Banach space and denote by $B(X)$ the algebra of
all bounded linear operators on $X$. The spectrum in any Banach algebra
is denoted by $\si(.)$. In $B(X)$, the spectrum has several
important subsets. In what follows $\si_p(A)$ denotes the point
spectrum of the operator $A\in B(X)$, that is,
\[
\si_p(A)=
\{ \lambda \in \C \, : \, A-\lambda I \text{ is noninjective} \}
\]
and $\si_s(A)$ denotes the surjectivity spectrum of $A$, that is,
\[
\si_s(A)=
\{ \lambda \in \C \, : \, A-\lambda I \text{ is nonsurjective} \}.
\]
If $x\in X$ and $f\in X^*$ ($X^*$ is the dual space of $X$), then $x\ot
f$ stands for the operator of rank at most one defined by
\[
(x\ot f)(y)=f(y)x \qquad (y\in X).
\]
Clearly, every finite rank operator $A \in B(X)$ is a finite linear
combination of such operators. On the finite rank elements of $B(X)$ one
can define the trace functional $\tr$ by
\[
\tr A=\sum_n f_n(x_n),
\]
where $A=\sum_n x_n \ot f_n$. Then $\tr$ is a well-defined linear
functional with
the property that $\tr(TA)=\tr(AT)$ for every finite rank operator $A\in
B(X)$ and for any $T\in B(X)$. For a matrix $A\in M_n(\C)$,
$A^t$ denotes the transpose of $A$.

If $X$ is a compact Hausdorff space, then let $C(X)$ denote the algebra
of all continuous complex valued functions on $X$. In this algebra the
spectrum of an element equals its range. If $f\in C(X)$, then $\supp f$
stands for the support of $f$, that is, $\supp f=\overline{\{ x\in X\,
:\, f(x)\neq 0\}}$.

Turning to our results and their proofs, we remark that on operator
algebras, besides linear maps preserving the spectrum one can also
consider such transformations which preserve some parts of the spectrum.
The following two results are of that type. In fact, they were motivated
by \cite[Theorem 3 and Theorem 4]{Semrl}, respectively.

\begin{theorem}\label{T:sielso}
Let $X$ be a Banach space and
let $\phi: B(X) \to B(X)$ be a surjective function with the property
that
\begin{equation}\label{E:alap2}
\si_p(\phi(A)\phi(B))=\si_p(AB) \qquad (A,B\in B(X)).
\end{equation}
If $X$ is infinite dimensional,
then there is an invertible linear operator $T\in B(X)$ such that
either
\[
\phi(A)=TAT^{-1} \qquad (A\in B(X))
\]
or
\[
\phi(A)=-TAT^{-1} \qquad (A\in B(X)).
\]
If $\phi: M_n(\C) \to M_n(\C)$ is a surjective function
satisfying \eqref{E:alap2}, then we have the following possibilities:
\begin{itemize}
\item[(a)]
there is an invertible matrix $T_1\in M_n(\C)$ such that
\[
\phi(A)=T_1 A T_1^{-1} \qquad (A\in M_n(\C));
\]
\item[(b)]
there is an invertible matrix $T_2\in M_n(\C)$ such that
\[
\phi(A)=-T_2 A T_2^{-1} \qquad (A\in M_n(\C));
\]
\item[(c)]
there is an invertible matrix $T_3\in M_n(\C)$ such that
\[
\phi(A)=T_3 A^t T_3^{-1} \qquad (A\in M_n(\C));
\]
\item[(d)]
there is an invertible matrix $T_4\in M_n(\C)$ such that
\[
\phi(A)=-T_4 A^t T_4^{-1} \qquad (A\in M_n(\C)).
\]
\end{itemize}
\end{theorem}

\begin{proof}
We first show that $\phi$ is injective. Indeed, if $A,A' \in B(X)$ are
such that $\phi(A)=\phi(A')$, then from \eqref{E:alap2} we obtain
that $\si_p(AB)=\si_p({A'}B)$ for every $B\in B(X)$.
This implies that
\begin{equation}\label{E:sielso}
\si_p( Ax \ot f)=\si_p(A'x \ot f) \qquad (x\in X, f\in X^*).
\end{equation}
It is an easy fact that if $\dim X\geq 2$, then
\begin{equation}\label{E:simasodik}
\si_p(y \ot g)=\{ 0, g(y) \} \qquad (y\in X, g\in X^*).
\end{equation}
Since in the one-dimensional case our statement is trivial, in what
follows we assume that $\dim X\geq 2$.
From \eqref{E:sielso} we infer that
$f(Ax)=f(A'x)$ for every $x\in X, f\in X^*$. It
follows that $A=A'$ which proves the injectivity of $\phi$.

Observe that $\phi$ preserves the rank-one operators.
In fact, this follows from the following characterization of rank-one
elements of $B(X)$.
The operator $A\in B(X)$ has rank one if and only if
$A\neq 0$, $0\in \si_p(TA)$ and $\# \si_p(TA) \leq 2$ for every $T\in
B(X)$ ($\#$ denotes cardinality). Observe that if $A\neq 0$, then
$\phi(A)\neq 0$.

Our next step is to show that $\phi$ is linear. The easiest way
to verify this is the use of the trace functional as follows.
Since the trace of a rank-one operator $x \ot f$ is $f(x)$,
we obtain from \eqref{E:alap2} and \eqref{E:simasodik} that
\begin{equation}\label{E:siharmadik}
\tr \phi(A)\phi(B)=\tr AB
\end{equation}
for every $A\in B(X)$ and rank-one operator $B\in B(X)$.
If $A,A'\in B(X)$ are arbitrary and $B\in B(X)$ is any rank-one
operator, then we compute
\[
\tr ((\phi(A)+\phi(A'))\phi(B))=
\tr \phi(A)\phi(B)+\tr \phi(A')\phi(B)=
\]
\[
\tr AB+\tr A'B=
\tr (A+A')B=
\tr \phi(A+A')\phi(B).
\]
By the arbitrariness of $B$ we obtain that $\phi$ is additive. One can
check that $\phi$ is homogeneous in a similar way.

So, $\phi$ is a linear bijection of $B(X)$ preserving the rank-one
operators.  The form of such transformations is well-known. It follows
from
the argument in \cite{Hou} leading to \cite[Lemma 1.2]{Hou} that we have
two possibilities:
\begin{itemize}
\item[(i)] there exist bijective linear operators $T:X\to X$ and $S: X^*
\to X^*$ such that
\[
\phi(x\ot f)=Tx \ot Sf \qquad (x\in X, f\in X^*);
\]
\item[(ii)] there exist bijective linear operators $T':X^*\to X$ and
$S': X \to X^*$ such that
\[
\phi(x\ot f)=T'f \ot S'x \qquad (x\in X, f\in X^*).
\]
\end{itemize}
Suppose first that we have (i). According to \eqref{E:siharmadik}
we obtain
\[
(Sf)(Ty)\cdot (Sg)(Tx)=f(y)g(x) \qquad (x,y \in X, f,g \in X^*).
\]
Consequently, there is a scalar $\lambda$ such that
\[
(Sg)(Tx)=\lambda g(x) \qquad (x \in X, g \in X^*).
\]
By the closed graph theorem we readily obtain that the
bijective linear operators
$T,S$ are bounded and hence we infer that $T^*S=\lambda I\in B(X^*)$.
Thus,
$S=\lambda {T^*}^{-1}=\lambda {T^{-1}}^*$ and this implies that
$\phi(A)=\lambda TAT^{-1}$ for every finite rank operator $A\in B(X)$.
Using the property \eqref{E:alap2} of $\phi$ we
have $\lambda^2=1$, that is, either $\lambda =1$ or $\lambda =-1$.
Suppose that $\lambda =1$.
Let $A\in B(X)$ be arbitrary. Pick any rank-one operator
$\phi(B)\in B(X)$. From \eqref{E:siharmadik} it follows that
\[
\tr \phi(A)\phi(B)=
\tr AB=
\tr(TAT^{-1}TBT^{-1})=
\tr(TAT^{-1}\phi(B)).
\]
By the arbitrariness of $\phi(B)$ we obtain that $\phi(A)=TAT^{-1}$ for
every $A\in B(X)$.

Assume now that we have (ii). Similarly to the case (i) one can prove
that ${T'}: X^*\to X$ is a bounded invertible linear operator and
\[
\phi(x\ot f)=\lambda {T'}(f\ot \tau(x)){T'}^{-1} \qquad (x\in X, f\in X^*)
\]
where $\tau$ denotes the natural embedding of $X$ into $X^{**}$. Since
$(x\ot f)^*=f\ot \tau(x)$, we obtain that in this case $\phi$ is of the
form $\phi(A)=\lambda {T'}A^*{T'}^{-1}$ for every finite rank operator $A\in B(X)$.
Just as above, we infer that $\lambda =\pm 1$ and then obtain the form
of $\phi$ on the whole $B(X)$. To see that in the infinite dimensional
case this second possibility (ii) cannot occur, we refer to
\cite[Theorem 3]{Semrl} stating that on an infinite
dimensional Banach space $X$ every point
spectrum preserving surjective linear map is an automorphism (not an
antiautomorphism). Since, as it can be seen, $\phi$ or $-\phi$ satisfies
these conditions, that result applies.

To verify that the finite dimensional case is different, that is, (ii)
can really occur, we remark
that in that case the injectivity, surjectivity, bijectivity of an
operator are all equivalent and that
it is true for any elements $A,B$ of any
Banach algebra that $\si(AB)\setminus \{0\}=\si(BA)\setminus \{0\}$.
Consequently, for every $A,B \in M_n(\C)$ we have  $\si_p(A^t
B^t)=\si_p((BA)^t)=\si_p(BA)=\si_p(AB)$. The proof is complete.
\end{proof}

Considering the surjectivity spectrum we have a similar result which
follows.

\begin{theorem}\label{T:simasodik}
Let $H$ be an infinite dimensional Hilbert space and
let $\phi: B(H) \to B(H)$ be a surjective function with the property
that
\begin{equation}\label{E:alap3}
\si_s(\phi(A)\phi(B))=\si_s(AB) \qquad (A,B\in B(H)).
\end{equation}
Then there is an invertible linear operator $T\in B(H)$ such that
either
\[
\phi(A)=TAT^{-1} \qquad (A\in B(H))
\]
or
\[
\phi(A)=-TAT^{-1} \qquad (A\in B(H)).
\]
\end{theorem}

\begin{proof}
One can argue in a very similar way as in our first result. This can
be done since, by the Fredholm alternative, for any
finite rank operator (in fact, even for any compact operator) $A\in
B(X)$ we have
\[
\si_p(A)\setminus \{ 0\}=\si(A)\setminus \{ 0\}=\si_s(A)\setminus
\{0\}.
\]
To exculde the appearence of the second possibility (ii) in the
proof of Theorem~\ref{T:sielso} choose a nonsurjective isometry $V$ on
$H$. Let $A=V^*$ (the Banach space adjoint of $V$) and set $B=V$. Then
we see that $AB$ is invertible while $A^*B^*$ is not surjective.
So, $\si_s(A^*B^*)\neq \si_s(AB)$.
\end{proof}

Using the same argument once again, we have the following result.

\begin{theorem}\label{T:simasodika}
Let $H$ be an infinite dimensional Hilbert space and
let $\phi: B(H) \to B(H)$ be a surjective function with the property
that
\begin{equation}\label{E:alap3a}
\si(\phi(A)\phi(B))=\si(AB) \qquad (A,B\in B(H)).
\end{equation}
Then there is an invertible linear operator $T\in B(H)$ such that
either
\[
\phi(A)=TAT^{-1} \qquad (A\in B(H))
\]
or
\[
\phi(A)=-TAT^{-1} \qquad (A\in B(H)).
\]
\end{theorem}

Now we turn to a similar characterization of *-automorphisms.

\begin{theorem}\label{T:siharmadik}
Let $H$ be a Hilbert space and
let $\phi: B(H) \to B(H)$ be a surjective function with the property
that
\begin{equation}\label{E:alap}
\si(\phi(A)^*\phi(B))=\si(A^*B) \qquad (A,B\in B(H)).
\end{equation}
If $H$ is infinite dimensional,
then there are unitary operators $U,V \in B(H)$ such that
$\phi$ is of the form
\[
\phi(A)=UAV \qquad (A\in B(H)).
\]
If $\phi: M_n(\C) \to M_n(\C)$ is a surjective funtion satisfying
\eqref{E:alap}, then there are unitaries $U,V \in M_n(\C)$ such that
$\phi$ is either of the form
\[
\phi(A)=UAV \qquad (A\in M_n(\C))
\]
or of the form
\[
\phi(A)=UA^t V \qquad (A\in M_n(\C)).
\]
\end{theorem}

\begin{proof}
The linearity of $\phi$ can be proved in the very similar way as
above. Since the norm and the spectral radius of a selfadjoint operator
in $B(H)$ are equal, it follows from $\si(\phi(A)^*\phi(A))=\si(A^*A)$
that $\| \phi(A)\|^2=\| A\|^2$ $(A,B\in B(H))$. Consequently,
$\phi$ is a surjective linear isometry of $B(H)$. The form of such
transformations is well-known. Namely, to every surjective linear
isometry $\psi$ there exist unitaries $U,V\in B(H)$ such that $\psi$ is
either of the form
\[
\psi(A)=UAV \qquad (A\in B(H))
\]
or of the form
\[
\psi(A)=UA^t V \qquad (A\in B(H)).
\]
If $H$ is of infinite dimension, then the appearence of
this second possibility can be excluded just as in the last part of the
proof of Theorem~\ref{T:simasodik}.
\end{proof}

We next treat our problems in the case of the function algebra $C(X)$ on
a compact Hausdorff space $X$.

\begin{theorem}\label{T:sinegyedik}
Let $X$ be a first countable compact Hausdorff space. If
$\phi :C(X) \to C(X)$ is a surjective function with the property that
\begin{equation}\label{E:alap5}
\si(\phi(f)\phi(g))=\si(fg) \qquad (f,g \in C(X)),
\end{equation}
then there exist a homeomorphism $\varphi :X \to X$
and a continuous function $\tau : X\to \{ -1, 1\}$
such that
\[
\phi(f)(x)=\tau(x) f(\varphi(x)) \qquad (x \in X, f\in C(X)).
\]
\end{theorem}

\begin{proof}
We have $\si(\phi(1)^2)=\si(1)$.
The spectrum of an element of $C(X)$ equals its range.
Therefore, $\phi(1)^2=1$ and
considering the transformation $f\mapsto \phi(1)\phi(f)$, we can and do
assume that our function $\phi$ satisfies $\phi(1)=1$.

We obtain from \eqref{E:alap5} that $\si(\phi(f))=\si(f)$ for every
$f\in C(X)$. So, $\phi$ preserves the range of functions. This implies
that $\phi$ maps real functions to real functions and it sends
nonnegative functions to nonnegative functions.

We prove that $\phi$ is injective. This will follow from
the following characterization of the equality between functions. Let
$f,g\in C(X)$. Then $f=g$ if and only if $\si(fh)=\si(gh)$ for every
nonnegative function $h\in C(X)$. Indeed, suppose that $f(x_0) \neq
g(x_0)$ for some $x_0\in X$. We can assume that $|f(x_0)|\leq
|g(x_0)|$. Let $D$ be an open disk centered at $f(x_0)$
which does not contain $g(x_0)$ and let $U$ be a neighbourhood of $x_0$
such that $f(x)\in D$ for every $x\in U$.  Let $h: X\to
[0,1]$ be a continuous function such that $\supp h\subset U$ and
$h(x_0)=1$. Such a function exists by Urysohn's lemma. Then we obtain
$\si (fh) \subset [0,1]D$
but $g(x_0)h(x_0)=g(x_0) \notin [0,1]D$. Therefore, $ \si(fh) \neq
\si(gh)$. So, we have the injectivity of $\phi$.
Therefore, $\phi$ and $\phi^{-1}$ are bijective functions having the
same properties concerning the spectrum.

Our next claim is that $\phi$ preserves the usual ordering between real
functions.
This will follow from the following characterization of that
ordering.
If $f,g \in C(X)$ are real functions, then $f\leq g$ if and only if
\begin{itemize}
\item[(a)]
$hg\leq c \Longrightarrow hf\leq c$ for every $0\leq h\in C(X)$ and
$c\in \R$

\hskip -1.3cm and

\item[(b)]
$hf\geq c \Longrightarrow hg\geq c$ for every $0\leq h\in C(X)$ and
$c\in \R$.
\end{itemize}
To see this, suppose that $f(x_0)>g(x_0)$ for some $x_0 \in X$. Clearly,
there exists a positive number $\epsilon$ such that either $f(x_0)$ does
not belong to
the $\epsilon$-neighbourhood of $[0,g(x_0)]$ or $g(x_0)$ does not belong
to the $\epsilon$-neighbourhood of $[0,f(x_0)]$.
Suppose that we have the first possibility. Choose a continuous function
$h:X \to [0,1]$ for which $h(x_0)=1$,  and the support of $h$ is a
subset of a neighbourhood of $x_0$ in which $g$ takes its
values in the $\epsilon$-neighbourhood of $g(x_0)$. Then
we find that $\si (hg)$ is a subset of the $\epsilon$-neighbourhood of
$[0, g(x_0)]$
but $\si (hf)$ is not a subset of that set. It is easy to see that there
is a constant $c$ such that $hg \leq c$ but $hf \nleq c$.
So, the above characterization really holds and then we get
that $f\leq g$ if and only if $\phi(f)\leq \phi(g)$.

Observe that by \eqref{E:alap5} we have $fg=0$ if and only if
$\phi(f)\phi(g)=0$. Now, if $f,g\geq 0$ and $fg=0$, then we find that
\[
\phi(f+g)=\phi(\max \{f, g\})=\max \{\phi(f), \phi(g)\}=
\phi(f)+\phi(g).
\]

To any point $x\in X$ there exists a continuous function $f_x :X \to
[0,1]$
such that $f_x(x)=1$ and $f_x(y) <1$ if $y\neq x$. In fact, by the first
countability of $X$ there is a sequence $\{ U_n\}$ of neighbourhoods of
$x$
which forms a base of neighbourhoods of that point. For every $n\in \N$
there is a continuous function
$f_n :X \to [0,1]$ whose support is in $U_n$ and $f_n(x)=1$. Now, set
$f_x=\sum_n \frac{1}{2^n} f_n$. This function fulfils our requirements.
Denote by $\mathcal F_x$ the set of all such function $f_x$. We assert
that if $f_x \in \mathcal F_x$, then $\phi(f_x)$ belongs to $\mathcal
F_{\varphi(x)}$ for some $\varphi(x)\in X$. In fact, since $\phi$
preserves
the range of functions, it follows that $\phi(f_x)$ maps into $[0,1]$
and it takes the value 1. Suppose that $\phi(f_x)$ equals 1 at two
different points, say $y,z$. It follows that there are functions $g'\in
\mathcal F_y$ and $h'\in \mathcal F_z$ such that $g'h'=0$ and $g'+h'\leq
\phi(f_x)$. Let $g=\phi^{-1}(g')$ and $h=\phi^{-1}(h')$. Then we have
$gh=0$ and by the previous sections of the
proof we infer that $g+h=\phi^{-1}(g'+h')\leq f_x$. Since $g+h$ takes
the
value 1 at at least two points, the same must be true for $f_x$ which is
a contradiction. This means that $\phi(f_x)\in \mathcal F_{\varphi(x)}$
for some $\varphi(x)\in X$.
We next show that the point $\varphi(x)$ does not depend on the
particular
choice of $f_x$. Indeed, let $f_x,f'_x \in \mathcal F_x$. Then
$\max \{f_x,f'_x\} \in \mathcal F_x$ and this implies that
$\max \{\phi(f_x),\phi(f'_x)\} \in \mathcal F_y$ for some $y\in X$. This
proves that $\phi(f_x)$ and $\phi(f'_x)$ take their maximum at the same
point.
So, we have a function $\varphi :X \to X$ such that
$f_x \in \mathcal F_x$ implies $\phi(f_x)\in \mathcal F_{\varphi(x)}$.
Since $\phi$ and $\phi^{-1}$ share the same properties, we obtain that
$\varphi$ is a bijection.

We assert that $\phi$ is homogeneous. Let $f\in C(X)$ and $\lambda \in
\C$. For any $0\leq h\in C(X)$ we have
\[
\si( \lambda \phi(f)\phi(h))=
\lambda \si(\phi(f)\phi(h))=
\lambda \si(fh)=
\si((\lambda f)h)=
\si(\phi(\lambda f)\phi(h))
\]
which implies that $\phi(\lambda f)=\lambda \phi(f)$.

Let $0\leq f\in C(X)$, $x\in X$ and let $f(x)=\lambda$. There exists
$f_x \in \mathcal F_x$ such that $\lambda f_x \leq f$. Then we have
$\lambda \phi(f_x)=\phi(\lambda f_x)\leq \phi(f)$. This gives us that
\[
f(x)=\lambda =\lambda \phi(f_x)(\varphi(x))\leq \phi(f)(\varphi(x)).
\]
Since $\phi^{-1}$ has the same properties as $\phi$, it follows that
\[
\phi^{-1}(\phi(f))(\varphi^{-1}(\varphi(x)))\geq \phi(f)(\varphi(x)),
\]
that is, we also have $f(x)\geq \phi(f)(\varphi(x))$. Therefore, we
obtain
$\phi(f)(\varphi(x))=f(x)$ for every $x\in X$ and $0\leq f\in C(X)$.

We show that $\varphi$ is a homeomorphism. We need only to show that
$\varphi$ is continuous. Let $x_n$ be a sequence in $X$ converging to
the point $x\in X$. Suppose that $\varphi_n(x)$ does not converge to
$\varphi(x)$. Then there is a neighbourhood $U$ of $\varphi(x)$ such
that $\varphi_n(x) \notin U$ for infinitely many indices. Let $h':X \to
[0,1]$ be a continuous function with support in $U$ such that
$h'(\varphi(x))=1$. Let $h\in C(X)$ be such that $\phi(h)=h'$. Then we
have
$h(x_n)=h'(\varphi(x_n))=0$ for infinitely many $n$'s and this
contradicts
$h(x_n) \to h(x)=1$. So, $\varphi$ is a homeomorphism of $X$ and
we have $\phi(f)=f\circ \varphi^{-1}$ for every nonnegative $f\in C(X)$.

Finally, for any $f\in C(X)$ and $0\leq h\in C(X)$ we compute
\[
\si(\phi(f) \cdot h\circ \varphi^{-1})=
\si(\phi(f) \phi(h))=
\si(f h)=
\si( f\circ \varphi^{-1} \cdot h\circ \varphi^{-1})
\]
which gives us that
$\phi(f)=f \circ \varphi^{-1}$.
This completes the proof.
\end{proof}

We turn to the second type of our preserver problems involving
involution. We have the following result.

\begin{theorem}\label{T:siotodik}
Let $X$ be a first countable compact Hausdorff space. If
$\phi :C(X) \to C(X)$ is a surjective function with the property that
\begin{equation}\label{E:alap6}
\si({\overline{\phi(f)}}\phi(g))=\si(\overline fg) \qquad (f,g \in
C(X)),
\end{equation}
then there exist a homeomorphism $\varphi: X\to X$ and a function $\tau
\in C(X)$ of modulus 1 such that
\[
\phi(f)(x)=\tau(x) f(\varphi(x)) \qquad (x\in X, f\in C(X)).
\]
\end{theorem}

\begin{proof}
Similarly as in the proof of Theorem~\ref{T:sinegyedik} one can verify
that $\phi$ is injective. Indeed, if $f,g\in C(X)$ are such that
$\phi(f)=\phi(g)$, then we have
\[
\si(\overline fh)=
\si({\overline{\phi(f)}}\phi(h))=
\si({\overline{\phi(g)}}\phi(h))=
\si(\overline gh)
\]
for every $h\in C(X)$ which implies
that $f=g$.

Observe that we have $|\phi(1)|^2=1$ which implies that
$\phi(1)$ is a function of modulus 1. Considering the transformation
$f \mapsto {\overline{\phi(1)}}\phi(f)$, we can and do assume that our
function $\phi$ satisfies $\phi(1)=1$. We have
\[
\si(\phi(g))=
\si({\overline{\phi(1)}}\phi(g))=\si({\overline{1}} g)=
\si(g)
\]
for every $g\in C(X)$. Therefore, $\phi$ is self-bijection of set
$C_\R(X)$ of all real valued continuous functions on $X$ which satisfies
\[
\si(\phi(f)\phi(g))=\si(fg) \qquad (f,g \in C_\R(X)).
\]
Since, as it turns out from the proof of the
previous result, Theorem~\ref{T:sinegyedik} remains valid for
the function algebra $C_\R(X)$ as well, we obtain that there is a
homeomorphism $\varphi: X\to X$ such that
\[
\phi(f)=f\circ \varphi \qquad (f\in C_\R(X)).
\]
If $f\in C(X)$, then we have
\[
\si(\phi(f) \cdot g\circ \varphi)=
\si(\phi(f) \phi(g))=
\si(f g)=
\si(f\circ \varphi \cdot g\circ \varphi)
\]
for every $g\in C_\R(X)$ which yields $\phi(f)=f\circ \varphi$.
The proof is complete.
\end{proof}

Finally, we present an application of our results. Let $\mathcal A$ be a
Banach algebra. The transformation $\phi :\mathcal A \to \mathcal A$ (no
linearity or continuity is assumed) is called a 2-local automorphism if
for every $x,y\in \mathcal A$ there exists an algebra
automorphism $\phi_{x,y}$ of $\mathcal A$ such that
$\phi(x)=\phi_{x,y}(x)$ and $\phi(y)=\phi_{x,y}(y)$.
Similarly, the transformation $\psi:\mathcal A \to \mathcal A$ is called
a 2-local isometry
if for every $x,y\in \mathcal A$ there exists a surjective
linear isometry $\psi_{x,y}$ of $\mathcal A$ such that
$\psi(x)=\psi_{x,y}(x)$ and $\psi(y)=\psi_{x,y}(y)$.
2-local maps were first studied by \v Semrl in \cite{Semrl2}.

Let $H$ be an infinite dimensional separable Hilbert space. It was
proved in \cite{Semrl2} that every 2-local automorphism of
$B(H)$ is an algebra automorphism of $B(H)$.
As for the function algebra $C(X)$ over a compact Hausdorff space $X$,
it follows from \cite[1.2. Theorem]{KoSlod} that every 2-local
automorphism of $C(X)$ is linear. Hence, applying our result
\cite[Theorem 2.2]{MoZal}, we see that if $X$ is a first countable
compact Hausdorff space, then every 2-local automorphism of $C(X)$ is an
algebra automorphism.

As for the isometry groups of the mentioned algebras, we refer to
\cite{Molnar} where we have proved that every 2-local isometry of any
$C^*$-subalgebra $\mathcal A$ of $B(H)$ which contains the ideal of
all compact operators and the identity operator
is linear. In particular, we obtained that every 2-local
isometry of $B(H)$ is a surjective linear isometry of $B(H)$.
Unfortunately, we do not know whether the analogue statement is true for
$C(X)$, $X$ being a first countable compact Hausdorff space. But it
follows from the form of the
surjective linear isometries of $C(X)$ given by Banach-Stone theorem and
from our result Theorem~\ref{T:siotodik} that every surjective 2-local
isometry is in fact a surjective linear isometry.

Referring back to \v Semrl's result on 2-local automorphisms of $B(H)$,
examining the proof of \cite[Theorem 1]{Semrl2}, it seems
essential that $H$
is a Hilbert space. It is a natural question that what can be stated for
Banach spaces. It follows from the form of the
automorphims of $B(X)$ (every algebra automorphism of $B(X)$ is
inner) and Theorem~\ref{T:sielso} that if $X$ is an
infinite dimensional Banach space and $\phi:B(X) \to B(X)$ is a
surjective 2-local
automorphism, then $\phi$ is an algebra automorphism of
$B(X)$. For an analogue result concerning linear (1-)local
automorphisms see \cite[Theorem 2.1]{LarSour}.

\bibliographystyle{amsplain}

\end{document}